\documentclass[a4paper,12pt]{article}
\usepackage{amsmath,amsfonts,amssymb,amsthm}
\newtheorem{theo}{Th\'eor\`eme}
\newtheorem*{coro}{Corollaire}
\newtheorem{lem}{Proposition}
\newtheorem*{ther}{Th\'eor\`eme}
\newcommand{\la}{\lambda}
\begin{document}
\title{Une $q$ - sp\'ecialisation pour les fonctions sym\'etriques monomiales}
\author{Michel Lassalle\\{\small Centre National de la Recherche Scientifique}\\
{\small Ecole Polytechnique}\\
{\small 91128 Palaiseau, France}\\
{\small e-mail: lassalle @ chercheur.com}}
\date{}
\maketitle

\begin{abstract}
We obtain the specialization of 
monomial symmetric functions on the alphabet $(a-b)/(1-q)$. This gives 
a remarkable algebraic identity, and four new developments for the Macdonald 
polynomial associated with a row. The proofs are given in the framework of 
$\lambda$-ring theory. 
\end{abstract}

\section{Introduction}

Dans l'\'etude des fonctions sym\'etriques, la th\'eorie des $\la$-anneaux 
est une m\'ethode particuli\`erement efficace, et pourtant peu 
utilis\'ee. On trouvera une illustration de cette th\'eorie
dans \cite{LL}. Le but de cet article est d'en pr\'esenter une 
nouvelle application.

Nous consid\'erons le probl\`eme suivant : si $f$ est une  
fonction sym\'etrique et $q$ une ind\'etermin\'ee, 
quelle est la valeur de $f(1,q,q^2,\ldots,q^{N-1})$~? 
Pour la plupart des fonctions sym\'etriques classiques
cette sp\'ecialisation est connue depuis tr\`es longtemps \cite{A,Li,Ma}. 

C'est le cas par exemple 
pour les fonctions de Schur, les sommes de puissances, les fonctions 
compl\`etes ou les fonctions \'el\'ementaires. 
Dans ces deux derniers cas cette sp\'ecialisation est 
classique : ce sont les polyn\^omes de Gauss.

Le but de cet article est de donner la sp\'ecialisation 
$f(1,q,q^2,\ldots,q^{N-1})$ lorsque $f$ est une fonction sym\'etrique 
monomiale. Ce r\'esultat n'\'etait pas encore connu. 
Plus g\'en\'eralement nous donnons la sp\'ecialisation des fonctions 
sym\'etriques monomiales sur l'alphabet $(a-b)/(1-q)$.

Il faut souligner que nous pouvons donner \emph{deux 
formulations} distinctes pour cette sp\'ecialisation. L'\'equivalence 
de ces deux expressions produit une identit\'e 
alg\'ebrique multivari\'ee qui est difficile \`a d\'emontrer directement. 
On a ainsi un nouvel exemple d'une situation o\`u la th\'eorie 
des $\la$-anneaux permet de d\'emontrer rapidement une identit\'e 
alg\'ebrique remarquable.

Nos deux r\'esultats et leurs d\'emonstrations s'expriment uniquement en termes de 
$\la$-anneaux. Mais ils poss\`edent des rapports \'etroits avec 
la th\'eorie des polyn\^omes de Macdonald  \cite{Ma}.

Soient $q$ et $t$ deux ind\'etermin\'ees, et consid\'erons l'alg\`ebre des 
fonctions sym\'etriques \`a coefficients rationnels en $q$ 
et $t$. Les polyn\^omes de Macdonald forment une base de cette alg\`ebre, 
index\'ee par les partitions. 

Notre r\'esultat principal permet d'obtenir quatre nouveaux
d\'eveloppements explicites pour le polyn\^ome de Macdonald
$P_{(n)}(q,t)$ associ\'e a une 
partition-ligne $(n)$. Pour cela nous introduisons deux bases 
naturelles de fonctions sym\'etriques ``d\'eform\'ees'', qui sont rest\'ees 
jusqu'ici peu \'etudi\'ees.

Enfin nous montrons que la sp\'ecialisation 
d'une fonction sym\'etrique monomiale sur l'alphabet $(1-t)/(1-q)$ 
est essentiellement un polyn\^ome en $q$ et $t$ \`a coefficients entiers positifs.
Il est possible que ce r\'esultat ``\`a la Macdonald'' ait d'int\'eressantes 
cons\'equences.

Donnons maintenant le plan de cet article. La Section 2 pr\'esente nos 
notations, et la Section 3 les \'el\'ements de th\'eorie des $\la$-anneaux 
dont nous aurons besoin. Ces sections sont presque int\'egralement reprises 
de \cite{LL}. 
La Section 4 \'enonce la premi\`ere formulation de notre r\'esultat 
principal, qui est d\'emontr\'ee \`a la Section 5. La Section 6 donne 
et d\'emontre la seconde formulation. 
La Section 7 met en \'evidence une identit\'e alg\'ebrique remarquable 
et pr\'esente quelques unes de ses cons\'equences.
La Section 8 explicite un polyn\^ome en $q$ et $t$
\`a coefficients entiers positifs.
La Section 9 introduit les polyn\^omes de 
Macdonald en mettant l'accent sur une pr\'esentation en termes de 
$\la$-anneaux.  
La Section 10 donne les nouveaux d\'eveloppements de $P_{(n)}(q,t)$ annonc\'es. 

L'auteur remercie Alain Lascoux pour son aide amicale.

\section{Notations}

Une partition $\la$ est une suite d\'ecroissante finie d'entiers positifs. On dit 
que le nombre $n$ d'entiers non nuls est la longueur de $\la$. On note
$\la  = ( \la_1,\ldots,\la_n)$ et $n = l(\la)$. On dit que 
$|\la| = \sum_{i = 1}^{n} \la_i$
est le poids de $\la$, et pour tout entier $i\geq1$ que 
$m_i(\la)  = \textrm{card} \{j: \la_{j}  = i\}$
est la multiplicit\'e de $i$ dans $\la$. On appelle part de $\la$ 
tout entier $i$ tel que $m_i(\la)\neq 0$. On pose 
\[z_\la  = \prod_{i \ge  1} i^{m_i(\lambda)} m_i(\lambda) !  .\]

Soit  $\mathbf{Sym}$  l'alg\`ebre des fonctions sym\'etriques. 
Nous choisissons les notations 
les plus r\'epandues, c'est-\`a-dire celles de \cite{Ma},
et non celles de \cite{LS}, bien que celles de \cite{LS} soient plus 
adapt\'ees aux $\lambda$-anneaux.

Soit $A=\{a_1,a_2,a_3,\ldots\}$ un ensemble de variables, qui peut \^etre infini 
(nous dirons que $A$ est un alphabet). On introduit les 
fonctions g\'en\'eratrices
\[ E_u(A) = \prod_{a\in A} (1 +ua) \quad , \quad 
 H_u(A)=  \prod_{a\in A}  \frac{1}{1-ua}
\quad , \quad P_u(A)= \sum_{a\in A} \frac{a}{1-ua} \]
dont le d\'eveloppement d\'efinit les fonctions sym\'etriques \'el\'ementaires
$e_{k}(A)$, les fonctions compl\`etes 
$h_{k}(A)$ et les sommes de puissances $p_{k}(A)$:

\[ E_u(A) =\sum_{k\geq0} u^k\, e_k(A)  \quad  , \quad 
 H_u(A) = \sum_{k\geq0} u^k\, h_k(A)
\quad  , \quad P_u(A)=\sum_{k\geq1} u^{k-1}p_k(A) .\]

Lorsque l'alphabet $A$ est infini, chacun de ces trois ensembles de fonctions forme une base alg\'ebrique de 
$\mathbf{Sym}[A]$, l'alg\`ebre des fonctions sym\'etriques sur $A$ (c'est-\`a-dire 
que ses \'el\'ements sont alg\'ebriquement ind\'ependants). 

On peut donc d\'efinir l'alg\`ebre $\mathbf{Sym}$ des fonctions sym\'etriques, 
sans r\'ef\'erence \`a l'alphabet $A$, comme l'alg\`ebre sur 
$\mathbf{Q}$ engendr\'ee par les fonctions $e_k$, $h_k$  
ou $p_k$.  

Pour toute partition $\mu$, on d\'efinit les fonctions $e_\mu$, 
$h_\mu$ ou $p_\mu$ en posant
\[f_{\mu}= \prod_{i=1}^{l(\mu)}f_{\mu_{i}}=\prod_{k\geq1}f_k^{m_{k}(\mu)},\]
o\`u $f_{i}$ d\'esigne respectivement $e_i$, $h_i$  ou $p_i$. 
Les fonctions $e_\mu$, $h_\mu$, $p_\mu$ forment une base 
lin\'eaire de l'alg\`ebre $\mathbf{Sym}$.

On a la formule de Cauchy
\[e_{n} = \sum_{\left|{\mu }\right| = n}(-1)^{n-l(\mu)} \frac{p_{\mu}}{z_{\mu}}\]
soit encore
\[h_{n} = \sum_{\left|{\mu }\right| = n} \frac{p_{\mu}}{z_{\mu}}.\]

Pour toute partition $\mu$, on peut d\'efinir les fonctions sym\'etriques monomiales 
$m_\mu$ et les fonctions de Schur $s_\mu$, qui forment \'egalement une base 
lin\'eaire de l'alg\`ebre $\mathbf{Sym}$. La fonction sym\'etrique monomiale 
$m_{\mu}$ est la somme de tous 
les mon\^omes diff\'erents ayant pour exposant une permutation de $\mu$.

Si $A$ et $B$ sont deux alphabets, on d\'efinit la somme $A+B$ 
et la diff\'erence $A-B$ de ces deux alphabets en posant
\begin{equation}
\begin{split}
H_u(A+B)=H_u(A)\,H_u(B) \quad&,\quad E_u(A+B)=E_u(A)\,E_u(B)\\
H_u(A-B)=H_u(A)\,{H_u(B)}^{-1} \quad&,\quad 
E_u(A-B)=E_u(A)\,{E_u(B)}^{-1}.
\end{split} 
\end{equation}

\section{$\la$-anneaux}

Nous allons utiliser le fait que 
l'anneau des polyn\^omes poss\`ede une structure de $\la$-anneau. 
Un $\la$-anneau est un anneau commutatif avec unit\'e muni 
d'op\'erateurs qui v\'erifient certains axiomes. Nous renvoyons le lecteur \`a  
\cite{K} pour la th\'eorie g\'en\'erale, et au chapitre 2  de \cite{Pr} 
pour son application \`a l'analyse multivari\'ee.

Nous n'utiliserons cette th\'eorie que dans le cadre \'el\'ementaire 
suivant. Soit $A=\{a_1,a_2,a_3,\ldots\}$ un alphabet quelconque. On consid\`ere 
l'anneau  $\mathbf{R}[A]$ des polyn\^omes en $A$ \`a coefficients r\'eels. 
La structure de  $\la$-anneau de $\mathbf{R}[A]$ consiste \`a 
d\'efinir une action de $\mathbf{Sym}$ sur $\mathbf{R}[A]$.

\subsection{Action de $\mathbf{Sym}$}

Les fonctions $p_{k}$ formant un syst\`eme de g\'en\'erateurs alg\'ebriques 
de $\mathbf{Sym}$, \'ecrivant tout polyn\^ome sous la forme $\sum_{c,U} 
c \, U$, avec $c$ constante r\'eelle et $U$ un mon\^ome en 
$(a_{1},a_{2},a_{3},\ldots)$, on d\'efinit une action de $\mathbf{Sym}$ 
sur $\mathbf{R}[A]$, not\'ee $[ \ . \ ]$, 
en posant 
\[ p_{k} [\sum_{c,U} c \, U]=\sum_{c,U} c \, U^{k}.\]

Pour tous polyn\^omes $P,Q \in \mathbf{R}[A]$ on en d\'eduit 
imm\'ediatement $p_{k} [PQ]=p_{k} [P] p_{k} [Q]$ et
$p_{\mu} [PQ]=p_{\mu} [P] p_{\mu} [Q]$.

L'action ainsi d\'efinie s'\'etend \`a tout \'el\'ement de $\mathbf{Sym}$. Ainsi on a
\[ E_u \,[\sum_{c,U} c\, U]= \prod_{c,U} (1+u\, U)^c \quad , \quad
H_u \,[\sum_{c,U} c \, U]= \prod_{c,U} (1-u\, U)^{-c},\]
et aussi 
\begin{equation}
h_{k} [P] =(-1)^k e_{k} [-P].
\end{equation}

On notera le comportement diff\'erent des 
constantes $c\in \mathbf{R}$ et des mon\^omes $U$ :
\begin{equation}
\begin{split}   
p_{k} [c] = c \quad , \quad  h_{k} [c] =\binom{c+k-1}{k} \quad , \quad  
e_{k} [c] = \binom{c}{k} \quad \\
p_k [U] = U^k = h_k [U] \quad , \quad  
e_{k} [U] = 0,\, i>1 \quad , \quad e_1 [U] = U.
\end{split}
\end{equation}

Il est plus correct de caract\'eriser les ``mon\^omes'' $U$ comme \'el\'ements {\it  
de rang 1} (i.e. les $U\neq 0,1$ tels que $e_{k}[U] = 0 \ \forall k>1$),
et les ``constantes'' $c\in \mathbf{R}$ comme les \'el\'ements invariants 
par les $p_k$ (on dira aussi \'el\'ement {\it de type binomial}).

Lorsqu'on utilise la th\'eorie des $\lambda$-anneaux pour d\'emontrer 
une identit\'e alg\'ebrique, il est donc \textit{toujours n\'ecessaire de 
pr\'eciser le statut de chaque \'el\'ement}. Dans cet article nous n'utiliserons 
que des \'el\'ements de rang 1.

\subsection{Extension aux s\'eries formelles}

On remarquera que si $a_{1},a_{2},\ldots,a_{N}$ sont des \'el\'ements de 
rang 1, alors  
\[p_{k} [a_{1}+a_{2}+\ldots+a_{N}]=a_{1}^k+a_{2}^k+\ldots+a_{N}^k\]
est la valeur de la somme de puissances $p_{k} 
(a_1,a_2,\ldots,a_{N})$. 

Pour tout alphabet 
$A=\{a_1,a_2,a_3,\ldots\}$, on note $A^\dag =\sum_{i}a_{i}$ la 
somme de ses \'el\'ements. Lorsque $A$ est form\'e d'\'el\'ements de rang 
1, on a ainsi pour toute fonction sym\'etrique $f$,
\begin{equation}
f[A^\dag]=f(A).
\end{equation}

En particulier si $q$ est de rang 1, on a 
\[p_{k} (1,q,q^2,q^3,\ldots,q^{N-1})= p_{k}\, [\sum_{i=0}^{N-1} q^i].\]
Il est naturel de vouloir \'ecrire
\[\sum_{i=0}^{N-1} q^i=\frac{1-q^N}{1-q},\] 
et d'\'etendre ainsi l'action de $\mathbf{Sym}$ aux fonctions rationnelles.  
Il est \'egalement naturel de consid\'erer un alphabet infini 
$(1,q,q^2,q^3,\ldots)$, de vouloir sommer la s\'erie 
\[\sum_{i\geq0}q^i=\frac{1}{1-q},\] 
et d'\'etendre ainsi l'action de $\mathbf{Sym}$ aux s\'eries formelles \`a 
coefficients r\'eels.

Pour cela on pose
\[ p_{k} \left( \frac{\sum c \,U}{\sum d \,V}\right) = 
  \frac{\sum c \,U^k}{\sum d\, V^k}  \ , \]
avec $c,d$ constantes r\'eelles et $U,V$ des mon\^omes en 
$(a_{1},a_{2},a_{3},\ldots)$.

L'action ainsi d\'efinie s'\'etend \`a tout \'el\'ement de $\mathbf{Sym}$. 
On munit ainsi l'anneau des s\'eries formelles \`a 
coefficients r\'eels d'une structure de $\lambda$-anneau.

\subsection{Formulaire}

Les relations fondamentales suivantes sont des 
cons\'equences directes des relations (3). Pour tous $P,Q$ on a d'abord
\begin{equation*}
\begin{split}
h_{n}[P+Q]&= \sum_{k=0}^{n} h_{n-k} [P] \,h_k [Q]\\
e_{n}[P+Q]&= \sum_{k=0}^{n} e_{n-k} [P] \,e_k [Q].
\end{split}
\end{equation*}
Soit de mani\`ere \'equivalente
\begin{equation}
\begin{split}
H_u[P+Q]=H_u[P]\,H_u[Q] \quad&,\quad E_u[P+Q]=E_u[P]\,E_u[Q]\\
H_u[P-Q]=H_u[P]\,{H_u[Q]}^{-1} \quad&,\quad 
E_u[P-Q]=E_u[P]\,{E_u[Q]}^{-1}.
\end{split} 
\end{equation}
Ces relations g\'en\'eralisent les d\'efinitions (1).

Si $P$ est de rang 1 et $Q$ arbitraire, on a
\[e_{n} [PQ] =P^n e_{n}[Q].\]
Si $P$ et $Q$ sont de rang 1, $PQ$ est donc de rang 1, et on a
\begin{equation}
E_u[PQ]=1+u\, PQ.
\end{equation}

Pour tous $P,Q$ 
on a les formules de Cauchy suivantes
\begin{equation}
\begin{split}
h_n [PQ]&= \sum_{\left|{\mu }\right| = n} 
\frac{1}{z_{\mu}}p_{\mu} [P] \, p_{\mu} [Q]\\
&=\sum_{\left|{\mu }\right| = n} m_{\mu} [P] \, h_{\mu} [Q]\\
&=\sum_{\left|{\mu }\right| = n} s_{\mu} [P] \, s_{\mu} [Q].
\end{split}
\end{equation}
Ou de mani\`ere \'equivalente :
\begin{equation}
\begin{split}
e_{n} [PQ]&= \sum_{\left|{\mu }\right| = n} 
\frac{(-1)^{n-l(\mu)}}{z_{\mu}} p_{\mu} [P] \, p_{\mu} [Q]\\
&=\sum_{\left|{\mu }\right| = n} m_{\mu} [P] \, e_{\mu} [Q]\\
&=\sum_{\left|{\mu }\right| = n} s_{\mu} [P] \, s_{\mu'} [Q],
\end{split}
\end{equation}
o\`u $\mu'$ d\'esigne la partition transpos\'ee de $\mu$.

\subsection{$q$-calcul}

Pour toute ind\'etermin\'ee $a$ on note
\[{(a;q)}_{n}=\prod_{i=0}^{n-1}(1-aq^i)\quad , \quad
{(a;q)}_{\infty}=\prod_{i\geq0}(1-aq^i)\]
qu'on consid\`ere comme s\'erie formelle en $a$ et $q$.

Soient trois \'el\'ements $a,b,q$ et un alphabet 
$X=\{x_1,x_2,\ldots,x_N\}$. On suppose tous ces Õel\'ements de rang 1. 
Les relation (5) et (6) impliquent
\begin{equation*}
\begin{split}
E_u \, \left[\frac{aX^\dag}{1-q}\right]&=\prod_{i\geq0} E_u \, [aq^iX^\dag] 
=\prod_{k=1}^{N} \, \prod_{i\geq0} E_u \, [aq^ix_k]\\
&=\prod_{k=1}^{N} \, \prod_{i\geq0} (1+uaq^ix_k)
=\prod_{k=1}^{N}{(-uax_k;q)}_{\infty}.
\end{split}
\end{equation*}
De m\^eme on a
\begin{equation*}
\begin{split}
H_u \,\left[\frac{aX^\dag}{1-q}\right]&=\prod_{i\geq0} H_u \, [aq^iX^\dag] 
=\prod_{k=1}^{N} \, \prod_{i\geq0} H_u \, [aq^ix_k]\\
&=\prod_{k=1}^{N} \, \prod_{i\geq0} \frac{1}{1-uaq^ix_k}
=\prod_{k=1}^{N} \frac{1}{{(uax_k;q)}_{\infty}}.
\end{split}
\end{equation*}

On en d\'eduit
\begin{equation*}
H_1 \,\left[\frac{a-b}{1-q}X^\dag\right]= 
H_1 \,\left[\frac{aX^\dag}{1-q}\right] 
{\left(H_1 \,\left[\frac{bX^\dag}{1-q}\right]\right)}^{-1} = \prod_{k=1}^{N} 
\frac{{(bx_k;q)}_{\infty}}{{(ax_k;q)}_{\infty}}.
\end{equation*}
Pour tout entier $n\ge 0$ on note d\'esormais
\begin{equation}
g_{n}(X;q,t)=h_{n}\left[\frac{1-t}{1-q}\,X^\dag\right].
\end{equation}
On a la s\'erie g\'en\'eratrice
\begin{equation}
H_1\left[\frac{1-t}{1-q}\,X^\dag\right]= \sum_{n\ge0} g_{n}(X;q,t)
= \prod_{i=1}^{N} 
\frac{{(tx_i;q)}_{\infty}}{{(x_i;q)}_{\infty}}.
\end{equation}
Les deux  propri\'et\'es suivantes sont des cons\'equences 
imm\'ediates de ce qui pr\'ec\`ede. Cependant nous en donnons 
une d\'emonstration directe \`a titre d'exemple. 

\begin{lem}
On a
\[\sum_{n\ge 0}q^n g_{n}(X;q,t) = H_1\,\left[q\,\frac{1-t}{1-q}\,X^\dag\right] = 
\left(\prod_{i=1}^N \frac{1-x_i}{1-tx_i}\right)\,
H_1\,\left[\frac{1-t}{1-q}\,X^\dag\right].\]
\end{lem}
\begin{proof}[Preuve]
On peut \'ecrire
\begin{equation*}
H_1\,\left[q\,\frac{1-t}{1-q}\,X^\dag\right] =
H_1 \,\left[\frac{1-t}{1-q}X^\dag+(t-1)X^\dag\right].
\end{equation*}
En appliquant (4) ceci devient
\[H_1\,\left[q\,\frac{1-t}{1-q}\,X^\dag\right]  =H_1[(t-1)X^\dag] \, H_1 
\,\left[\frac{1-t}{1-q}X^\dag\right].\]
Mais on a
\[H_1[(t-1)X^\dag] = H_1[tX^\dag]{H_1[X^\dag]}^{-1}
=\prod_{i=1}^N \frac{1-x_i}{1-tx_i}.\]
\end{proof}

\begin{lem}
On a 
\[H_1\left[\frac{q-t}{1-q}\,X^\dag\right]=
\prod_{i=1}^N (1-x_i)\,
H_1\left[\frac{1-t}{1-q}\,X^\dag\right].\]
\end{lem}
\begin{proof}[Preuve]
On peut \'ecrire
\begin{equation*}
H_1\,\left[\frac{q-t}{1-q}\,X^\dag\right]
=H_1 \,\left[\frac{1-t}{1-q}X^\dag-X^\dag\right]
\end{equation*}
En appliquant (4) ceci devient
\[H_1\left[\frac{q-t}{1-q}\,X^\dag\right]
= H_1 \,\left[\frac{1-t}{1-q}X^\dag\right] {H_1[X^\dag]}^{-1}.\]
\end{proof}

\section{Notre r\'esultat principal}

Etant donn\'ee une partition $\mu$, on note $C_{\mu}$ l'ensemble des 
multi-entiers distincts obtenus par 
permutation des parts de $\mu$. On dit \'egalement que $c \in C_{\mu}$ est 
un ``d\'erangement'' de $\mu$. Pour tout multi-entier 
$c=(c_1,\ldots,c_{l(\mu)})\in C_{\mu}$, on note $[c_i]=\sum_{k\le i} c_k$ la somme 
partielle d'ordre $i$.

Soient $a,b,q$ trois \'el\'ements de rang $1$. 
Nous consid\'erons l'alphabet $A$ tel que 
\[A^\dag = \frac{a-b}{1-q}.\]
L'alphabet $A$ est la diff\'erence, au sens de (1), des
deux alphabets infinis $\{a, aq, aq^2,\ldots\}$ et
$\{b, bq, bq^2,\ldots\}$.

Soit $m_{\mu}$ la fonction sym\'etrique 
monomiale associ\'ee \`a la partition $\mu$. 
Nous explicitons la valeur de $m_{\mu}[A^\dag]$.
En particulier pour $a=1$ et $b=q^N$, compte-tenu de (4),
notre r\'esultat donne la valeur de 
\[m_{\mu}\left[\frac{1-q^N}{1-q}\right] = m_{\mu}(1,q,\ldots,q^{N-1}).\]
Et pour $a=1$ et $b=0$ celle de
\[m_{\mu}\left[\frac{1}{1-q}\right] = m_{\mu}(1,q,q^2,q^3,\ldots).\]

\begin{theo}
Soient $a,b,q$ trois \'el\'ements de rang 1. Pour toute 
partition $\mu$ on a
\begin{equation}
m_{\mu}\left[\frac{a-b}{1-q}\right]=\sum_{c\in C_{\mu}} \prod_{i=1}^{l(\mu)}
\frac{\displaystyle{a^{c_i}q^{[c_{i-1}]}-b^{c_i}}}{\displaystyle{1-q^{[c_i]}}}.
\end{equation}
\end{theo}

Le corollaire suivant \'etait connu : voir \cite{A}, chapitre 2, et 
l'exemple 1.2.5 de \cite{Ma}. Compte-tenu de (10), il s'agit du ``th\'eor\`eme 
de Heine'' classique
\[H_1 \,\left[\frac{1-t}{1-q}\,x\right]= 
\frac{{(tx;q)}_{\infty}}{{(x;q)}_{\infty}}=
\sum_{n\ge 0} \frac{{(t;q)}_n}{{(q;q)}_n} \,x^n.\]

\begin{coro}
Pour tout entier $n$ on a
\begin{equation*}
\begin{split}
e_n\left[\frac{a-b}{1-q}\right]&=\prod_{i=1}^{n}
\frac{aq^{i-1}-b}{1-q^i}\\
h_n\left[\frac{a-b}{1-q}\right]&=\prod_{i=1}^{n}
\frac{a-bq^{i-1}}{1-q^i}.
\end{split}
\end{equation*}
\end{coro}
\begin{proof}[Preuve du corollaire]
La premi\`ere relation est la transcription du th\'eor\`eme pour la 
partition-colonne $\mu=1^n$. On a alors $m_{1^n}=e_{n}$ . Il n'y a qu'un seul 
d\'erangement de $\mu$, avec 
$c_{i}=1$ et $[c_i]=i$. La seconde relation s'en d\'eduit par (2), 
en \'echangeant $a$ et $b$.
\end{proof}

On remarquera que le Th\'eor\`eme 1 est v\'erifi\'e lorsque $\mu$ est une 
partition-ligne $(n)$. On a alors $m_{(n)}=p_{n}$. Il n'y a qu'un seul 
d\'erangement de $\mu$, avec $c_1=[c_1]=n$. 
Le Th\'eor\`eme 1 redonne dans ce cas la relation
\[p_n\left[\frac{a-b}{1-q}\right]=
\frac{a^n-b^n}{1-q^n},\]
ce qui est pr\'ecis\'ement la d\'efinition de l'action de $p_{n}$.

On note d\'esormais
\[Z_{\mu}(a,b,q) = \sum_{c\in C_{\mu}} \prod_{i=1}^{l(\mu)}
\frac{\displaystyle{a^{c_i}q^{[c_{i-1}]}-b^{c_i}}}{\displaystyle{1-q^{[c_i]}}}.\]
Pour toute part 
$i$ de $\mu$, on note $\mu \backslash \{i\}$ la partition de longueur $l(\mu)-1$ 
obtenue de $\mu$ par soustraction de $i$.
\begin{lem}
Pour toute partition $\mu$ on a
\begin{equation}
(1-q^{|\mu|})Z_{\mu}(a,b,q) = 
\sum_{\begin{subarray}{c}i\\m_{i}(\mu) \neq 0\end{subarray}} 
(a^iq^{|\mu|-i}-b^i)Z_{\mu \backslash \{i\}}(a,b,q).
\end{equation}
\end{lem}
\begin{proof}[Preuve]
On consid\`ere tous les d\'erangements de $\mu$ dont la derni\`ere 
composante est $c_{l(\mu)}=i$. On a alors $[c_{l(\mu)-1}]=|\mu|-i$ et
$[c_{l(\mu)}]=|\mu|$. Par construction la somme de toutes ces 
contributions est exactement
\[Z_{\mu \backslash \{i\}}(a,b,q) \frac 
{\displaystyle{(a^iq^{|\mu|-i}-b^i)}}{\displaystyle{1-q^{|\mu|}}}.\]
\end{proof}

Partant du cas initial $\mu=(n)$ la relation (12) d\'etermine uniquement
$Z_{\mu}$ par r\'ecurrence sur la longueur $l(\mu)$. Le Th\'eor\`eme 1
sera donc d\'emontr\'e si l'on \'etablit que le membre de gauche de 
(11) satisfait la m\^eme relation de r\'ecurrence.

Les deux membres de (11) \'etant clairement homog\`enes de degr\'e $|\mu|$,
il suffit de d\'emontrer le th\'eor\`eme dans le cas particulier 
$a=1$, ce que nous supposerons d\'esormais. 

\section{D\'emonstration}

Nous sommes ainsi conduits \`a d\'emontrer le Th\'eor\`eme 1 sous la forme 
suivante.
\begin{theo}
Soient $q$ et $t$ deux \'el\'ements de rang 1. Pour toute 
partition $\mu$ on a
\begin{equation*}
(1-q^{|\mu|}) \, m_{\mu}\left[\frac{1-t}{1-q}\right] = 
\sum_{\begin{subarray}{c}i\\m_{i}(\mu) \neq 0\end{subarray}} 
(q^{|\mu|-i}-t^i) \, m_{\mu \backslash \{i\}} \left[\frac{1-t}{1-q}\right].
\end{equation*}
\end{theo}

\begin{proof}[Preuve]
Soit $X=\{x_1,x_2,\ldots,x_N\}$ un alphabet de cardinal $N$ dont les 
\'el\'ements sont de rang 1. Compte-tenu de la relation (4),
de la d\'efinition (9) et de  la formule de Cauchy (7), on a imm\'ediatement
\[g_{n}(X;q,t) = h_n \left[\frac{1-t}{1-q} X^\dag\right]
= \sum_{\left|{\mu }\right| = n} m_{\mu} \left[\frac{1-t}{1-q}\right] \, 
h_{\mu} (X).\]

En identifiant les parties homog\`enes de chaque membre, 
le Th\'eor\`eme 2 est donc \'equivalent \`a la relation suivante
\begin{multline*}
\sum_{n\ge0}(1-q^n)\,g_{n}(X;q,t) =\\
\left(\sum_{r\ge1} h_r(X)\right)\left(\sum_{n\ge0}q^ng_{n}(X;q,t)\right)
-\left(\sum_{r\ge1} t^r h_r(X)\right) \left(\sum_{n\ge0}g_{n}(X;q,t)\right).
\end{multline*}
Soit encore
\begin{multline*}
\sum_{n\ge0}(1-q^n) \, g_{n}(X;q,t) =\\
\left(\prod_{i=1}^{N}  \frac{1}{1-x_i} -1\right)
\left(\sum_{n\ge0}q^ng_{n}(X;q,t)\right)
-\left(\prod_{i=1}^{N}  \frac{1}{1-tx_i} -1\right) 
\left(\sum_{n\ge0}g_{n}(X;q,t)\right).
\end{multline*}

On applique alors la Proposition 1. Le membre de gauche
peut s'\'ecrire 
\[\left(1-\prod_{i=1}^N \frac{1-x_i}{1-tx_i}\right)\,
H_1\left[\frac{1-t}{1-q}\,X^\dag\right].\]
Et le membre de droite s'\'ecrit
\[\left[\left(\prod_{i=1}^{N}  \frac{1}{1-x_i} -1\right)
\left(\prod_{i=1}^N \frac{1-x_i}{1-tx_i}\right)
-\left(\prod_{i=1}^{N}  \frac{1}{1-tx_i}-1\right)\right]\, 
H_1\left[\frac{1-t}{1-q}\,X^\dag\right].\]
D'o\`u l'assertion.
\end{proof}

\section{Seconde formulation}

Il est remarquable que nous puissions donner une seconde formulation de notre 
r\'esultat principal.

\begin{theo}
Soient $a,b,q$ trois \'el\'ements de rang 1. Pour toute partition $\mu$ on a  
\[m_{\mu}\left[\frac{a-b}{1-q}\right]=
\sum_{c\in C_{\mu}} \prod_{i=1}^{l(\mu)}
\frac{\displaystyle{a^{c_i}q^{(l(\mu)-i)c_i}-b^{c_i}}}{\displaystyle{1-q^{[c_i]}}}
.\]
\end{theo}

Notre d\'emonstration du Th\'eor\`eme 3 est exactement parall\`ele 
\`a celle du Th\'eor\`eme 1. On note
\[W_{\mu}(a,b,q) = \sum_{c\in C_{\mu}} \prod_{i=1}^{l(\mu)}
\frac{\displaystyle{a^{c_i}q^{(l(\mu)-i)c_i}-b^{c_i}}}{\displaystyle{1-q^{[c_i]}}}
.\]

\begin{lem}
Pour toute partition $\mu$ on a
\begin{equation*}
(1-q^{|\mu|})W_{\mu}(a,b,q) = 
\sum_{\begin{subarray}{c}i\\m_{i}(\mu) \neq 0\end{subarray}} 
(a^i-b^i)W_{\mu \backslash \{i\}}(qa,b,q).
\end{equation*}
\end{lem}
\begin{proof}[Preuve]
On consid\`ere tous les d\'erangements de $\mu$ dont la derni\`ere 
composante est $c_{l(\mu)}=i$. On a alors
$[c_{l(\mu)}]=|\mu|$. Par construction la somme de toutes ces 
contributions est exactement
\[W_{\mu \backslash \{i\}}(qa,b,q) \, 
\frac {\displaystyle{(a^i-b^i)}}{\displaystyle{1-q^{|\mu|}}}.\]
\end{proof}

Partant du cas initial \'evident $\mu=(n)$, la Proposition 4 d\'etermine 
uniquement $W_{\mu}$ par r\'ecurrence sur la longueur $l(\mu)$. Le Th\'eor\`eme 3
sera donc d\'emontr\'e si l'on \'etablit que
$m_{\mu}[(a-b)/(1-q)]$ satisfait la m\^eme relation de r\'ecurrence.
Par homog\'eneit\'e il suffit de le d\'emontrer dans le cas particulier 
$a=1$, ce que nous supposerons d\'esormais. 

Nous pouvons donc d\'emontrer le Th\'eor\`eme 3 sous la forme 
suivante.
\begin{theo}
Soient $q$ et $t$ deux \'el\'ements de rang 1. Pour toute 
partition $\mu$ on a
\begin{equation*}
(1-q^{|\mu|}) \, m_{\mu}\left[\frac{1-t}{1-q}\right] = 
\sum_{\begin{subarray}{c}i\\m_{i}(\mu) \neq 0\end{subarray}} 
(1-t^i) \, m_{\mu \backslash \{i\}} \left[\frac{q-t}{1-q}\right].
\end{equation*}
\end{theo}

\begin{proof}[Preuve]
Soit $X=\{x_1,x_2,\ldots,x_N\}$ un alphabet de cardinal $N$ dont les 
\'el\'ements sont de rang 1. Compte-tenu de (4) la 
formule de Cauchy (7) s'\'ecrit
\[h_n \left[\frac{u-t}{1-q} X^\dag\right]
= \sum_{\left|{\mu }\right| = n} m_{\mu} \left[\frac{u-t}{1-q}\right] \, 
h_{\mu} (X).\]
On va choisir $u=1$ et $u=q$.

En identifiant les parties homog\`enes de chaque membre, 
le Th\'eor\`eme 3 est \'equivalent \`a la relation suivante
\begin{multline*}
\sum_{n\ge0}(1-q^n)\,h_n \left[\frac{1-t}{1-q} X^\dag\right] =\\
\left(\sum_{r\ge1} h_r(X)-\sum_{r\ge1} t^r h_r(X)\right)
\left(\sum_{n\ge0}h_n \left[\frac{q-t}{1-q} X^\dag\right]\right).
\end{multline*}
Soit encore
\begin{multline*}
\sum_{n\ge0}(1-q^n) \, g_{n}(X;q,t) =\\
\left(\prod_{i=1}^{N}  \frac{1}{1-x_i} -
\prod_{i=1}^{N}  \frac{1}{1-tx_i}\right) 
\left(\sum_{n\ge0}h_n \left[\frac{q-t}{1-q} X^\dag\right]\right).
\end{multline*}
Par la Proposition 1 le membre de gauche peut s'\'ecrire
\[\left(1-\prod_{i=1}^N \frac{1-x_i}{1-tx_i}\right)\,
H_1\left[\frac{1-t}{1-q}\,X^\dag\right] .\]
Et par la Proposition 2 le membre de droite peut s'\'ecrire
\[\left(\prod_{i=1}^{N}  \frac{1}{1-x_i} -
\prod_{i=1}^{N}  \frac{1}{1-tx_i}\right) 
\prod_{i=1}^N (1-x_i)\,
H_1\left[\frac{1-t}{1-q}\,X^\dag\right].\]
D'o\`u l'assertion.
\end{proof}

\section{Une identit\'e remarquable}

La comparaison des Th\'eor\`emes 1 et 3 produit  
l'identit\'e remarquable suivante.
\begin{theo}
Pour toute partition $\mu$ on a
\begin{equation*}
\sum_{c\in C_{\mu}} \prod_{i=1}^{l(\mu)}
\frac{\displaystyle{a^{c_i}q^{[c_{i-1}]}-b^{c_i}}}{\displaystyle{1-q^{[c_i]}}}
=\sum_{c\in C_{\mu}} \prod_{i=1}^{l(\mu)}
\frac{\displaystyle{a^{c_i}q^{(l(\mu)-i)c_i}-b^{c_i}}}{\displaystyle{1-q^{[c_i]}}}.
\end{equation*}
\end{theo}

Un cas particulier int\'eressant est obtenu en faisant $a=q$ et $b=1$.
\begin{lem}
Pour toute partition $\mu$ on a
\[\sum_{c\in C_{\mu}} \prod_{i=1}^{l(\mu)}
\frac{\displaystyle{1-q^{(l(\mu)-i+1)c_i}}}{\displaystyle{1-q^{[c_i]}}}
=\frac{l(\mu)!}{\prod_{i}m_{i}(\mu)!} .\]
\end{lem}

On peut prendre la limite de ce r\'esultat lorsque $q$ tend vers $1$.
\begin{lem}
Pour toute partition $\mu$ on a
\[\frac{1}{z_{\mu}} =\sum_{c\in C_{\mu}} \, \prod_{i=1}^{l(\mu)} 
\frac{1}{\displaystyle{[c_i]}}.\]
\end{lem}
Nous retrouvons ainsi un \'enonc\'e de  Littlewood (\cite{Li}, p. 85) qui l'a 
d\'emontr\'e par r\'ecurrence.

Soient $X=\{x_1,x_2,\ldots,x_n\}$ et
$Y=\{y_1,y_2,\ldots,y_n\}$ deux alphabets de cardinal 
$n$. On note $S_n$ le groupe des permutations de $n$ lettres. 
Le groupe $S_n$ op\`ere sur les fonctions rationnelles 
en $X$ et $Y$ par l'action diagonale
\[f^{\sigma}(x_1,\ldots,x_n,y_1,\ldots,y_n) =
f(x_{\sigma(1)},\ldots,x_{\sigma(n)},
y_{\sigma(1)},\ldots,y_{\sigma(n)}).\]

Par homog\'en\'eit\'e et parce que les ind\'etermin\'ees  
$x_i=q^{\mu_i}$, $y_i={(bq/a)}^{\mu_i}$ 
sont ind\'ependantes, l'\'egalit\'e du Th\'eor\`eme 5 est en fait
\emph{\'equivalente} \`a l'identit\'e multivari\'ee suivante,
qui est une propri\'et\'e des fonctions rationnelles. 
\begin{theo}
On a
\begin{multline*}
\sum_{\sigma \in S_n} 
\left(\frac{y_1-x_1}{1-x_1}\,\frac{y_2-x_1x_2}{1-x_1x_2}\,
\frac{y_3-x_1x_2x_3}{1-x_1x_2x_3}\cdots
\frac{y_n-x_1x_2\cdots x_n}{1-x_1x_2\cdots x_n}\right)^\sigma =\\
\sum_{\sigma \in S_n}
\left(\frac{y_1-{x_1}^{n}}{1-x_1}\,\frac{y_2-{x_2}^{n-1}}{1-x_1x_2}\,
\frac{y_3-{x_3}^{n-2}}{1-x_1x_2x_3}\cdots
\frac{y_n-x_n}{1-x_1x_2\cdots x_n}\right)^\sigma.
\end{multline*}
\end{theo}

Le caract\`ere remarquable de cette identit\'e est d\'ej\`a apparent 
sur le cas particulier $Y=(1,1,\ldots,1)$. 
\begin{lem}
Pour tout alphabet $X=\{x_1,x_2,\ldots,x_n\}$ on a
\[\sum_{\sigma \in S_n} 
\left(\frac{1-{x_1}^{n}}{1-x_1}\,\frac{1-{x_2}^{n-1}}{1-x_1x_2}\,
\frac{1-{x_3}^{n-2}}{1-x_1x_2x_3}\cdots
\frac{1-{x_n}}{1-x_1x_2\cdots x_n}\right)^\sigma = n!.\]
\end{lem}
On en d\'eduit la propri\'et\'e suivante.
\begin{lem}
Pour tout alphabet $X=\{x_1,x_2,\ldots,x_n\}$ on a
\[\sum_{\sigma \in S_n} 
\left(\frac{1}{x_1(x_1+x_2)\cdots(x_1+x_2+\cdots + x_n)}\right)^\sigma =
\prod_{i=1}^n \frac{1}{x_i}.\]
\end{lem}
\begin{proof}[Preuve]
On consid\`ere l'identit\'e de la Proposition 7, dans laquelle on 
substitute  $q^{x_i}$ \`a $x_i$. Elle devient
\[\sum_{\sigma \in S_n} 
\left(\frac{1-{q}^{nx_1}}{1-q^{x_1}}\,\frac{1-{q}^{(n-1)x_2}}{1-q^{x_1+x_2}}\,
\frac{1-{q}^{(n-2)x_3}}{1-q^{x_1+x_2+x_3}}\cdots
\frac{1-{q}^{x_n}}{1-q^{x_1+x_2\cdots+x_n}}\right)^\sigma= n!.\]
On obtient l'\'enonc\'e en prenant la limite $q\rightarrow 1$.
\end{proof}

Alain Lascoux a obtenu une preuve directe de cette identit\'e 
multivari\'ee, en utilisant les diff\'erences divis\'ees. 
Il nous a \'egalement montr\'e que le Th\'eor\`eme 6 
\'enonce l'\'egalit\'e de deux statistiques sur le groupe des 
permutations. Nous pr\'esentons maintenant l'essentiel de ses remarques.

Etant donn\'ee une permutation de $n$ lettres 
$\sigma\in S_n$, soit $\Gamma(\sigma)$ l'ensemble de ses cycles. Pour tout cycle 
$\gamma= (\gamma_1,\ldots,\gamma_k) \subset\{1,2,\ldots,n\}$, notons
$|\gamma|=\sum_{i=1}^k \gamma_i$.
Alors on sait (\cite{LS}, \S 1.2.7) que pour toute partition $\mu$ on a
\[(\prod_{i\ge1}m_i(\mu)!) \ m_\mu =\sum_{\sigma\in S_{l(\mu)}}
(-1)^{l(\mu)-\textrm{card}(\Gamma(\sigma))} 
\prod_{\gamma \in \Gamma(\sigma)}p_{|\gamma|}.\]
Par exemple on a $m_{kl}=p_kp_l-p_{k+l}$, chacun des termes correspondant 
aux deux cycles $\{k\},\{l\}$ de $\{k,l\}$ et au cycle $\{k,l\}$ de $\{l,k\}$.
On en d\'eduit imm\'ediatement
\[(\prod_{i\ge1}m_i(\mu)!) \ m_\mu \left[\frac{a-b}{1-q}\right]
=\sum_{\sigma\in S_{l(\mu)}}
(-1)^{l(\mu)-\textrm{card}(\Gamma(\sigma))} 
\prod_{\gamma\in \Gamma(\sigma)}
\frac{a^{|\gamma|}-b^{|\gamma|}}{1-q^{|\gamma|}}.\]

Par homog\'en\'eit\'e et parce que les ind\'etermin\'ees  
$x_i=q^{\mu_i}$, $y_i={(bq/a)}^{\mu_i}$ 
sont ind\'ependantes, ceci implique imm\'ediatement le r\'esultat 
suivant.
\begin{theo}
Soient $X=\{x_1,x_2,\ldots,x_n\}$ et
$Y=\{y_1,y_2,\ldots,y_n\}$ deux alphabets de cardinal 
$n$. On a
\begin{multline*}
\sum_{\sigma \in S_n} 
\left(\frac{y_1-x_1}{1-x_1}\,\frac{y_2-x_1x_2}{1-x_1x_2}\,
\frac{y_3-x_1x_2x_3}{1-x_1x_2x_3}\cdots
\frac{y_n-x_1x_2\cdots x_n}{1-x_1x_2\cdots x_n}\right)^\sigma =\\
\sum_{\sigma \in S_n}
\prod_{\begin{subarray}{c} \gamma \in \Gamma(\sigma)\\
\gamma=(\gamma_1,\ldots,\gamma_k)\end{subarray}}
\frac{y_{\gamma_1}y_{\gamma_2}\cdots y_{\gamma_k}
-x_{\gamma_1}x_{\gamma_2}\cdots x_{\gamma_k}}
{1-x_{\gamma_1}x_{\gamma_2}\cdots x_{\gamma_k}}.
\end{multline*}
\end{theo}
Tandis que le Th\'eor\`eme 6 concerne la sym\'etrisation de deux fonctions 
rationnelles, le Th\'eor\`eme 7  fait intervenir la structure des cycles 
d'une permutation, ce qui est une information non imm\'ediate sur cette permutation. 

Nous donnons en Appendice une preuve directe du Th\'eor\`eme 7, par r\'ecurrence 
sur l'entier $n$. Nous ne connaissons pas de preuve directe du Th\'eor\`eme 6.

\emph{Exemple : } Le cas $n=1$ est trivial. Dans le cas $n=2$, les deux identit\'es des 
Th\'eor\`emes 6 et 7 s'\'ecrivent
\begin{multline*}
\frac{y_1-x_1}{1-x_1}\,\frac{y_2-x_1x_2}{1-x_1x_2}\,
+\frac{y_2-x_2}{1-x_2}\,\frac{y_1-x_1x_2}{1-x_1x_2}=\\
\frac{y_1-{x_1}^{2}}{1-x_1}\,\frac{y_2-{x_2}}{1-x_1x_2}\,
+\frac{y_2-{x_2}^{2}}{1-x_2}\,\frac{y_1-{x_1}}{1-x_1x_2}=\\   
\frac{y_1-x_1}{1-x_1}\, \frac{y_2-x_2}{1-x_2}\,
+\frac{y_1y_2-x_1x_2}{1-x_1x_2}.
\end{multline*}

\section{Coefficients entiers positifs}

Nous allons voir que dans l'\'enonc\'e du Th\'eor\`eme 3,
\[m_{\mu}\left[\frac{1-t}{1-q}\right]=
=\sum_{c\in C_{\mu}} \prod_{i=1}^{l(\mu)}
\frac{\displaystyle{q^{(l(\mu)-i)c_i}-t^{c_i}}}{\displaystyle{1-q^{[c_i]}}}\]
le membre de droite met en \'evidence un polyn\^ome en $q$ et 
$t$ \`a coefficients entiers positifs.

\begin{theo}
Pour toute partition $\mu$ on a 
\[m_{\mu}\left[\frac{1-t}{1-q}\right]=\frac{l(\mu)!}{\prod_{i}m_{i}(\mu)!}\,
\left(\prod_{i=1}^{l(\mu)} \frac{q^{i-1}-t}{1-q^i}\right) 
\frac{H_{\mu}(q,t)}{q^{|\mu|-l(\mu)}H_{\mu}\,(q,1/q)}.\]
o\`u  $H_{\mu}(q,t)$ est un polyn\^ome en $q$ et $t$, et 
$q^{|\mu|-l(\mu)}H_{\mu}(q,1/q)$ un polyn\^ome en $q$. Les coefficients de
$H_{\mu}(q,t)$ sont entiers positifs.
\end{theo}

\begin{proof}[Preuve] 
Consid\'erons le polyn\^ome $P_{\mu}(q)$ d\'efini par
\[P_{\mu}(q)=\prod_{k=1}^{l(\mu)}\, 
\prod_{1 \le i_1 < i_2 < \ldots < i_k \le l(\mu)}
\frac{\displaystyle{1-q^{\sum \mu_{i_{j}}}}}{1-q}.\]
C'est \'evidemment un polyn\^ome
en $q$ \`a coefficients entiers positifs. Posons
\[H_{\mu}(q,t)=P_{\mu}(q)
\sum_{c\in C_{\mu}} \prod_{i=1}^{l(\mu)}
\frac{\displaystyle{q^{(l(\mu)-i)c_i}-t^{c_i}}}
{\displaystyle{q^{l(\mu)-i}-t}}
\,\frac{1-q}{\displaystyle{1-q^{[c_i]}}}.\]
Il est clair qu'on d\'efinit ainsi un polyn\^ome
en $q$ et $t$ \`a coefficients entiers positifs.
On a imm\'ediatement
\[q^{|\mu|-l(\mu)}H_{\mu}(q,1/q)=P_{\mu}(q)
\sum_{c\in C_{\mu}} \prod_{i=1}^{l(\mu)}
\frac{\displaystyle{1-q^{(l(\mu)-i+1)c_i}}}
{\displaystyle{1-q^{l(\mu)-i+1}}}
\,\frac{1-q}{\displaystyle{1-q^{[c_i]}}}.\]
On en d\'eduit que $q^{|\mu|-l(\mu)}H_{\mu}(q,1/q)$ est un polyn\^ome
en $q$. D'autre part on a
\[H_{\mu}(q,t)=P_{\mu}(q) \, 
\left(\prod_{i=1}^{l(\mu)}\frac{1-q}{q^{i-1}-t}\right)
\sum_{c\in C_{\mu}} \prod_{i=1}^{l(\mu)}
\frac{\displaystyle{q^{(l(\mu)-i)c_i}-t^{c_i}}}
{\displaystyle{1-q^{[c_i]}}}.\]
Mais en appliquant la Proposition 5 on a aussi
\[\frac{l(\mu)!}{\prod_{i}m_{i}(\mu)!} \,P_{\mu}(q)=
q^{|\mu|-l(\mu)} \,
\left(\prod_{i=1}^{l(\mu)}\frac{1-q^i}{1-q}\right) \,
H_{\mu}(q,1/q).\]
D'o\`u l'\'enonc\'e.
\end{proof}

\emph{Exemple : } Dans le cas particulier d'une partition $\mu=(n,k)$, avec deux 
parts distinctes $n\neq k$, on a
\[P_{n,k}(q)=\frac{1-q^n}{1-q}\, \frac{1-q^k}{1-q} \,\frac{1-q^{n+k}}{1-q}.\]
Le polyn\^ome $H_{n,k}(q,t)$ est donn\'e par
\[H_{n,k}(q,t)=
\frac{q^n-t^n}{q-t}\,\frac{1-t^k}{1-t}\,\frac{1-q^k}{1-q}
+\frac{q^k-t^k}{q-t}\,\frac{1-t^n}{1-t}\,\frac{1-q^n}{1-q}.\]

Dans le cas g\'en\'eral, il serait int\'eressant de disposer d'une interpr\'etation 
combinatoire de $H_{\mu}(q,t)$.

\section{Polyn\^omes de Macdonald}

La r\'ef\'erence pour les polyn\^omes de Macdonald est le Chapitre 6 
de \cite{Ma}. Nous rappelons seulement ici les \'el\'ements dont nous aurons 
besoin, en mettant l'accent sur une pr\'esentation en termes de 
$\la$-anneaux.

Soient deux \'el\'ements $q,t$ et un alphabet 
$X=\{x_1,x_2,\ldots,x_N\}$. On suppose tous ces \'el\'ements de rang 1.
Pour tout $1\le i \le N$, on pose
\[A_i(X;t) =\prod_{\begin{subarray}{c}j=1\\j\neq i\end{subarray}}^N
\frac{tx_i-x_j}{x_i-x_j}.\]

On note $T_{x_i}$ l'op\'erateur de $q$-d\'eformation d\'efini par
\[T_{x_i}f (x_1,\ldots,x_N)= f(x_1,\ldots,qx_i,\ldots,x_N).\]
Les polyn\^omes de Macdonald $P_{\la}(X;q,t)$ 
sont les vecteurs propres de l'op\'erateur aux diff\'erences
\[D(X;q,t)= \sum_{i=1}^{N} \, A_i(X;t) \,T_{x_i}.\]
On a
\[D(X;q,t) \, P_{\la}(X;q,t)=  \left(\sum_{i=1}^N q^{\la_i}\, 
t^{N-i}\right) P_{\la}(X;q,t).\]

On peut munir l'alg\`ebre des fonctions sym\'etriques \`a coefficients 
rationnels en $q$ et $t$ d'un produit scalaire $<\, , \,>_{q,t}$ d\'efini par
\[<p_\la,p_\mu>_{q,t}=\delta_{\la \mu}\,z_{\la}\,
p_{\la} \left[\frac{1-q}{1-t}\right].\]
Les polyn\^omes de Macdonald $P_{\la}(X;q,t)$ forment une base 
orthogonale pour ce produit scalaire. 
Si on note $Q_{\la}(X;q,t)$ la base duale on a
\begin{equation*}
\begin{split}   
H_1\left[\frac{1-t}{1-q}\,X^\dag \,Y^\dag\right]&=
\sum_{\la}P_{\la}(X;q,t)\,Q_{\la}(Y;q,t)\\
&=\sum_{\la}h_\la\left[\frac{1-t}{1-q}\,X^\dag\right]\,m_{\la}(Y)\\
&=\sum_{\la} s_{\la} [(1-t)\,X^\dag]
\, s_{\la} \left[\frac{Y^\dag}{1-q}\right].
\end{split} 
\end{equation*}
o\`u les deux derni\`eres relations r\'esultent de la formule de Cauchy (7).

On sait (\cite{Ma}, relation (4.9), p. 323) que le polyn\^ome de Macdonald 
$P_{(n)}(X;q,t)$ est proportionnel \`a $g_{n}(X;q,t)$. Cependant dans \cite{Ma} 
ce r\'esultat n'est pas d\'emontr\'e directement. Il nous parait 
int\'eressant d'en pr\'esenter une d\'emonstration directe 
dans le cadre des $\la$-anneaux.
\begin{theo}
On a
\[D(X;q,t) \,g_{n}(X;q,t)= \left( q^n t^{N-1} + \frac{1-t^{N-1}}{1-t}\right) 
g_{n}(X;q,t).\]
\end{theo}
\begin{proof}[Preuve]
Nous donnons une preuve \'el\'ementaire, mais il s'agit d'un cas 
particulier du Th\'eor\`eme 2.1 de \cite{LaMo}, qui est beaucoup plus 
g\'en\'eral. Compte-tenu de la d\'efinition (9), il faut prouver
\[D(X;q,t) \,H_1\left[\frac{1-t}{1-q}\,X^\dag\right]=
t^{N-1}\sum_{n\ge 0}q^n g_{n}(X;q,t) +
\frac{1-t^{N-1}}{1-t} \, H_1\left[\frac{1-t}{1-q}\,X^\dag\right].\]
Compte-tenu de la Proposition 1, ceci est \'equivalent \`a
\[D(X;q,t) \,H_1\left[\frac{1-t}{1-q}\,X^\dag\right]
=\left(t^{N-1}\prod_{i=1}^N \frac{1-x_i}{1-tx_i}+
\frac{1-t^{N-1}}{1-t}\right) H_1\left[\frac{1-t}{1-q}\,X^\dag\right].\]

Mais on voit facilement que
\[T_{x_i}\,h_{n}\left[\frac{1-t}{1-q}\,X^\dag\right]=
h_{n}\left[\frac{1-t}{1-q}\,(X^\dag+(q-1)x_i)\right]=
h_{n}\left[\frac{1-t}{1-q}\,X^\dag+(t-1)x_i\right].\]
En appliquant (5) ceci s'\'ecrit
\begin{equation*}
\begin{split}
T_{x_i}\,H_1\left[\frac{1-t}{1-q}\,X^\dag\right]&=
H_1\left[\frac{1-t}{1-q}\,X^\dag+(t-1)x_i\right]\\
&=H_1\left[\frac{1-t}{1-q}\,X^\dag\right]H_1[(t-1)x_i]\\
&=H_1\left[\frac{1-t}{1-q}\,X^\dag\right]\frac{1-x_i}{1-tx_i}.
\end{split}
\end{equation*}
L'assertion est alors une cons\'equence imm\'ediate de la proposition suivante.
\end{proof}
\newpage
\begin{lem}
On a
\[\sum_{i=1}^{N} A_{i}(X;t) = \frac{1-t^N}{1-t}\]
\[\sum_{i=1}^{N} \frac{x_i}{1-tx_i} A_{i}(X;t) = 
\frac{t^{N-1}}{1-t} \left(1- \prod_{i=1}^{N} \frac{1-x_i}{1-tx_i}\right).\]
\end{lem}
\begin{proof}[Preuve]
Le principe est celui donn\'e dans l'exemple 6.3.2 (a) de \cite{Ma}.
Il suffit de choisir $u=0$ et $u=1/t$ dans l'identit\'e de
d\'ecomposition en \'el\'ements simples suivante
\[\prod_{i=1}^{N} \frac{tu-x_i}{u-x_i}=
(t-1)\sum_{i=1}^{N}\frac{x_i\,A_i(X;t)}{u-x_i} +t^N.\]
Cette relation est une interpolation de Lagrange.
D\'efinissons le r\'esultant de deux alphabets $A$ et $B$ par 
\[R(A,B)= \prod_{a\in A, b\in B} (a-b).\]
On rappelle  \cite{L} que si $f(a)$ est un polyn\^ome  ayant $a^N$ comme terme de 
plus haut degr\'e, on a
\[\sum_{a\in A}\frac{f(a)}{R(a,A - a)}=1\]
pour tout alphabet $A$ de cardinal $N+1$.
La relation pr\'ec\'edente n'est autre que cette identit\'e \'ecrite pour 
$A=X+u$ et $f(a) = R(a,X/t)$, c'est-\`a-dire
\[\sum_{i=1}^{N} \frac{x_i-x_i/t}{x_i-u}
\prod_{\begin{subarray}{c}j=1\\j\neq i\end{subarray}}^N
\frac{x_i-x_j/t}{x_i-x_j} +
\prod_{i=1}^N
\frac{u-x_i/t}{u-x_i}=1.\]
\end{proof}

\section{D\'eveloppements}

Les Th\'eor\`emes 1 et 3 permettent d'\'ecrire plusieurs d\'eveloppements 
explicites pour le polyn\^ome de Macdonald $g_{n}(X;q,t)$. A chaque fois, il s'agit 
d'une application \'el\'ementaire de la relation (4), de la d\'efinition 
(9) et des formules de Cauchy (7--8).

\subsection{Bases classiques}

Nous redonnons d'abord deux r\'esultats connus. Le premier est l'exemple 
6.8.8(a) de \cite{Ma}. On a
\begin{equation*}
\begin{split}
g_{n}(X;q,t)&= \sum_{|\mu| = n} 
\frac{1}{z_{\mu}}p_{\mu} \left[\frac{1-t}{1-q}\right] \, p_{\mu} [X^\dag]\\
&=\sum_{|\mu| = n} \frac{1}{z_{\mu}}\,\prod_{i=1}^{l(\mu)} 
\frac{1-t^{\mu_i}}{1-q^{\mu_i}} \, p_{\mu} (X).
\end{split}
\end{equation*}
Le second est l'exemple 6.2.1 de \cite{Ma}. On a
\begin{equation*}
\begin{split}
g_n(X;q,t)&= \sum_{|\mu| = n} 
h_\mu \left[\frac{1-t}{1-q}\right] \, m_\mu [X^\dag]\\
&= \sum_{|\mu| = n} 
\prod_{i=1}^{l(\mu)} \, 
\frac{{(t;q)}_{\mu_i}}{{(q;q)}_{\mu_i}} \, m_\mu (X)
\end{split}
\end{equation*}
o\`u la derni\`ere \'egalit\'e r\'esulte du Corollaire du Th\'eor\`eme 1.

Les relations suivantes sont nouvelles. On a
\begin{equation*}
\begin{split}
g_{n}(X;q,t)&= \sum_{|\mu| = n} 
m_\mu \left[\frac{1-t}{1-q}\right] \, h_\mu [X^\dag]\\
&=\sum_{|\mu| = n} Z_\mu (1,t,q) \, h_\mu (X).
\end{split}
\end{equation*}
Et de m\^eme
\begin{equation*}
\begin{split}
g_{n}(X;q,t)&= (-1)^n \sum_{|\mu| = n} 
m_\mu \left[\frac{t-1}{1-q}\right] \, e_\mu [X^\dag]\\
&= (-1)^n \sum_{|\mu| = n}  Z_\mu (t,1,q) \, e_\mu (X).
\end{split}
\end{equation*}

\subsection{Bases ``d\'eform\'ees''}

Pour toute partition $\mu$ on pose
\[E_{\mu}(X;t) = e_{\mu} [(1-t)X^\dag] \quad , \quad
H_{\mu}(X;t) = h_{\mu} [(1-t)X^\dag].\]

En appliquant les formules de Cauchy (7--8), nous obtenons le 
d\'eveloppement explicite de $g_{n}(X;q,t)$ sur ces bases
\begin{equation*}
\begin{split}
g_{n}(X;q,t)&= \sum_{|\mu| = n} m_{\mu} \left[\frac{1}{1-q}\right] \, 
h_{\mu} [(1-t)X^\dag]\\
&=\sum_{|\mu| = n} Z_{\mu} (1,0,q) \, H_{\mu}(X;t)
\end{split}
\end{equation*}
\begin{equation*}
\begin{split}
g_{n}(X;q,t)&=(-1)^n \sum_{|\mu| = n}  m_{\mu} \left[\frac{1}{q-1}\right] \, 
e_{\mu} [(1-t)X^\dag]\\
&=(-1)^n \sum_{|\mu| = n}  Z_{\mu} (0,1,q) \, E_{\mu}(X;t).
\end{split}
\end{equation*}

Consid\'erons l'endomorphisme $\omega_{q,t} : f \rightarrow \omega_{q,t}(f)$ 
d\'efini sur toute fonction sym\'etrique homog\`ene par
\[\omega_{q,t}(f)[X] = (-1)^{\mathrm{deg}(f)} f \left[\frac{q-1}{1-t}X\right].\]
Compte-tenu de (2) on a imm\'ediatement
\[\omega_{q,t}(g_{n}(X;q,t))=e_{n}(X)\]
\[\omega_{q,t}(E_{\mu}(X;t))= H_{\mu}(X;q)\quad,\quad 
\omega_{q,t}(H_{\mu}(X;t))= E_{\mu}(X;q).\]

\subsection{Formulaire}

Les fonctions $E_{n}(X;t)$ et $H_{n}(X;t)$, et donc les 
bases $E_{\mu}(X;t)$ et $H_{\mu}(X;t)$, sont explicitement 
connues. 

\begin{lem}
Pour tout entier $n \ge 1$ on a
\begin{equation*}
\begin{split}
E_{n}(X;t)&= {(-t)}^n g_{n}(X;0,1/t) = 
{(-1)}^n t^{n-N}(t-1) \sum_{i=1}^N A_{i}(X;t) x_i^n\\
H_{n}(X;t)&= g_{n}(X;0,t) = t^{N-1}(1-t) \sum_{i=1}^N A_{i}(X;1/t) x_i^n.
\end{split}
\end{equation*}
\end{lem}
\begin{proof}[Preuve]
Les premi\`eres \'egalit\'es sont \'evidentes. Les secondes 
r\'esultent de l'exemple 6.3.2 (a) de \cite{Ma}, qui se d\'emontre comme 
la Proposition 9.
\end{proof}

On en d\'eduit le d\'eveloppement des fonctions $E_{n}(X;t)$ ou $H_{n}(X;t)$ 
sur les bases classiques. En effet les formules de Cauchy de la Section 
10.1 impliquent imm\'ediatement
\begin{equation*}
\begin{split}
E_{n}(X;t)&= \sum_{|\mu| = n} 
\frac{(-1)^{n-l(\mu)}}{z_{\mu}} \prod_{i=1}^{l(\mu)} (1-t^{\mu_i}) \, 
p_{\mu} (X)\\
&=(-1)^n \sum_{|\mu| = n} Z_{\mu} (t,1,0) \, h_{\mu} (X)\\
&=\sum_{|\mu| = n} Z_{\mu} (1,t,0) \, e_{\mu} (X)\\
&= \sum_{|\mu| = n} {(-t)}^{n-l(\mu)}{(1-t)}^{l(\mu)} \, m_\mu (X).
\end{split}
\end{equation*}

Et de m\^eme
\begin{equation*}
\begin{split}
H_{n}(X;t)&=\sum_{|\mu| = n} 
\frac{1}{z_{\mu}} \prod_{i=1}^{l(\mu)} (1-t^{\mu_i}) \, p_{\mu} (X)\\
&=\sum_{|\mu| = n} Z_{\mu} (1,t,0) \, h_{\mu} (X)\\
&= (-1)^n \sum_{|\mu| = n} Z_{\mu} (t,1,0) \, e_{\mu} (X)\\
&= \sum_{|\mu| = n} {(1-t)}^{l(\mu)} \, m_\mu (X).
\end{split}
\end{equation*}

Inversement on a
\begin{equation*}
\begin{split}
h_{n}(X)&= g_{n}(X;q,q)\\
&=\sum_{|\mu| = n} Z_{\mu} (1,0,q) \, H_{\mu} (X;q)\\
&= (-1)^n \sum_{|\mu| = n} Z_{\mu} (0,1,q) \, E_{\mu} (X;q).
\end{split}
\end{equation*}

Et de m\^eme
\begin{equation*}
\begin{split}
e_{n}(X)&=\sum_{|\mu| = n} Z_{\mu} (1,0,q) \, E_{\mu} (X;q)\\
&= (-1)^n \sum_{|\mu| = n} Z_{\mu} (0,1,q) \, H_{\mu} (X;q).
\end{split}
\end{equation*}

\section{Appendice}

Nous donnons ici une preuve directe du Th\'eor\`eme 7.
\begin{ther}
Soient $X=\{x_1,x_2,\ldots,x_n\}$ et
$Y=\{y_1,y_2,\ldots,y_n\}$ deux alphabets de cardinal 
$n$. On a
\begin{multline*}
\sum_{\sigma \in S_n} 
\left(\frac{y_1-x_1}{1-x_1}\,\frac{y_2-x_1x_2}{1-x_1x_2}\,
\frac{y_3-x_1x_2x_3}{1-x_1x_2x_3}\cdots
\frac{y_n-x_1x_2\cdots x_n}{1-x_1x_2\cdots x_n}\right)^\sigma =\\
\sum_{\sigma \in S_n}
\prod_{\begin{subarray}{c} \gamma \in \Gamma(\sigma)\\
\gamma=(\gamma_1,\ldots,\gamma_k)\end{subarray}}
\frac{y_{\gamma_1}y_{\gamma_2}\cdots y_{\gamma_k}
-x_{\gamma_1}x_{\gamma_2}\cdots x_{\gamma_k}}
{1-x_{\gamma_1}x_{\gamma_2}\cdots x_{\gamma_k}}.
\end{multline*}
\end{ther}
\begin{proof}[Preuve]
Les deux membres de l'identit\'e sont lin\'eaires en $y_n$. Il suffit 
donc de la d\'emontrer pour $y_n=1$ et $y_n=x_n$. 

Soit $L_n$ (resp. $R_n$) le membre de gauche (resp. 
de droite). Par r\'ecurrence sur l'entier $n$, il suffit de d\'emontrer que pour 
$f_n=L_n$ et $f_n=R_n$, on a les deux relations 
\begin{multline}
f_n(x_1,\ldots,x_n;y_1,\ldots,y_{n-1},x_n)=\\
\sum_{i=1}^{n-1} f_{n-1}(x_1,\ldots,x_ix_n,\ldots,x_{n-1};
y_1,\ldots,y_ix_n,\ldots,y_{n-1})
\end{multline}
\begin{multline}
f_n(x_1,\ldots,x_n;y_1,\ldots,y_{n-1},1)=
f_{n-1}(x_1,\ldots,x_{n-1};y_1,\ldots,y_{n-1})+\\
\sum_{i=1}^{n-1} f_{n-1}(x_1,\ldots,x_ix_n,\ldots,x_{n-1};
y_1,\ldots,y_i,\ldots,y_{n-1}).
\end{multline}

A titre d'exemple, nous montrons (13) pour $L_n$ et (14) pour $R_n$. 
La v\'erification de (13) pour $R_n$ et (14) pour $L_n$ est
exactement identique. Nous la laissons au lecteur.

Pour $y_n=x_n$ seules les permutations avec $\sigma(1) \neq n$ contribuent 
au membre de gauche. Supposons qu'on a $\sigma(i) = n$ avec 
$i \neq 1$. Le terme
\[\frac{y_{\sigma(i-1)}-x_{\sigma(1)}\cdots x_{\sigma(i-1)}}
{1-x_{\sigma(1)}\cdots x_{\sigma(i-1)}}\,
\frac{y_{\sigma(i)}-x_{\sigma(1)}\cdots x_{\sigma(i)}}
{1-x_{\sigma(1)}\cdots x_{\sigma(i)}}\]
devient
\[\frac{y_{\sigma(i-1)}-x_{\sigma(1)}\cdots x_{\sigma(i-1)}}
{1-x_{\sigma(1)}\cdots x_{\sigma(i-1)}}\,
\frac{x_{n}(1-x_{\sigma(1)}\cdots x_{\sigma(i-1)})}
{1-x_{\sigma(1)}\cdots x_{\sigma(i)}} =
\frac{y_{\sigma(i-1)}x_n-x_{\sigma(1)}\cdots x_{\sigma(i-1)}x_n}
{1-x_{\sigma(1)}\cdots x_{\sigma(i-1)}x_n}.\]
Ce qui prouve (13) pour le membre de gauche.

Pour $y_n=1$ toutes les permutations contenant le cycle $(n)$ contribuent 
au terme $R_{n-1}(x_1,\ldots,x_{n-1};y_1,\ldots,y_{n-1})$. Toutes les 
autres permutations ont un cycle $\gamma=(\gamma_1,\ldots,\gamma_k)$ de la forme 
$(\delta n)$, avec $\delta=(\gamma_1,\ldots,\gamma_{k-1})$. La contribution de 
ce cycle est
\[\frac{y_{\gamma_1}y_{\gamma_2}\cdots y_{\gamma_k}
-x_{\gamma_1}x_{\gamma_2}\cdots x_{\gamma_k}}
{1-x_{\gamma_1}x_{\gamma_2}\cdots x_{\gamma_k}}=
\frac{y_{\gamma_1}y_{\gamma_2}\cdots y_{\gamma_{k-1}}
-x_{\gamma_1}x_{\gamma_2}\cdots x_{\gamma_{k-1}}x_n}
{1-x_{\gamma_1}x_{\gamma_2}\cdots x_{\gamma_{k-1}}x_n}.\]
Ce qui d\'emontre (14) pour le membre de droite.
\end{proof}

\end{document}